\documentclass[11pt]{article}
\usepackage{amsmath,amsthm,amsfonts,amssymb,mathrsfs,epsfig,bm}
\usepackage{color}
\usepackage{xspace}
\parskip 	\bigskipamount
\newtheorem{theorem}{Theorem}

\theoremstyle{remark}

\def\N{\mathbb{N}}

\def\R{\mathbb{R}}

\renewcommand{\phi}{\varphi}
\renewcommand{\epsilon}{\varepsilon}

\newcommand{\comp}{\raisebox{0.1ex}{\scriptsize $\circ$}}

\newcommand{\var}{\operatorname{var}}

\newcommand{\keywords}[1]{ \noindent {\footnotesize
             {\small \em Keywords and phrases.} {\sc #1} } }
\newcommand{\ams}[2]{  \noindent {\footnotesize
             {\small \em AMS {\rm 2000} subject classifications.
             {\rm Primary {\sc #1}; secondary {\sc #2}} } } }
\usepackage{amsmath,amsthm,amsfonts,amssymb,mathrsfs,epsfig,bm}

\begin{document}

\title{\bf A note on the convergence of renewal and regenerative processes
to a Brownian bridge
}
\author{
{\sc Serguei Foss}
\and
{\sc Takis Konstantopoulos} 
}
\date{12 August 2007}
\maketitle

\begin{abstract}
The standard functional central limit theorem for a renewal
process with finite mean and variance, results in a Brownian
motion limit.
This note shows how to obtain a Brownian bridge process
by a direct procedure that does not involve conditioning.
Several examples are also considered.

\vspace*{2mm}

\vspace*{2mm}
\ams{60F17}{60K05,60G50}
\end{abstract}

\section{The basic theorem}
In proving convergence results for a stochastic ordered graph on the
integers \cite{FK03}, we noticed that one can obtain a Donsker-like
theorem for Brownian bridge in a somewhat non-standard manner.
The result appears to be new. As it may be of potential
interest in some related areas (statistics, large deviations), 
we summarise it in this short note.

Consider a (possibly delayed) renewal process on $[0, \infty)$
with renewal epochs 
\[
0 < R_1 < R_2 < \cdots.
\] 
We assume that
$\{R_{n+1}-R_n\}_{n\ge 1}$ are i.i.d.\ with mean $\mu$
and variance $\sigma^2$, both finite.
Let 
\[
A_t:= \#\{n \ge 1: ~ R_n \le t\}
\] 
be the associated counting process.
The standard functional central limit theorem for a renewal
process, see, e.g., \cite{ASM03},
states that the sequence of processes $\xi_1, \xi_2, \ldots$, where
\[
\xi_n(t) := \frac{A_{nt} - \mu^{-1} nt}{\sqrt{n}}, 
\quad t \ge 0,
\]
converges weakly, as $n \to \infty$, to $\mu^{-3/2}\sigma W$,
where $W$ is a standard Brownian motion on $[0, \infty)$.
Weak convergence (denoted by $\Rightarrow$ below)
means weak convergence of probability measures
on the space $D[0,\infty)$ of functions which are right
continuous with left limits, equipped with the usual Skorokhod
topology (see, e.g., \cite{BILL68}, \cite{WHITT02}).

A standard Brownian bridge \cite[p.\ 84]{BILL68} $W^0$ is defined, 
in distribution, as a
standard Brownian motion $W$ on $[0,1]$, conditional on $W_1=0$, i.e.\
as the weak limit of the sequence of probability measures
\[
P(W \in \cdot \mid 0 \le W_1 \le 1/n), \quad n \in \N,
\]
as $n \to \infty$.
Often, when Brownian bridge is obtained as a limit
by a functional central limit theorem, there is an {\em explicit} underlying
conditioning that takes place. One first proves convergence
to a Brownian motion and uses conditioning to prove
convergence to a Brownian bridge.
Brownian bridges appear in limits of urn processes,
and also in limits of empirical distributions
\cite[Thm.\ 13.1]{BILL68}.

In this note we remark that it is possible
to obtain a Brownian bridge from a renewal process, without
the use of conditioning.

\begin{theorem}
\label{BBthm}
Define, for $u >0$,
\[
\eta_u(t) := \frac{R_{[t A_u]} - t u}{\sqrt{u}},
\quad 0 \le t \le 1.
\]
Considering $\eta_u$ as a random element of $D[0,1]$ (equipped with the 
topology of uniform convergence on compacta),
we have 
\[
\eta_u \Rightarrow \mu^{-1/2} \sigma W^0,
\quad  \text{ as $u\to\infty$},
\] 
where $W^0$ is a standard Brownian bridge.
\end{theorem}

Here, $[x]$ denotes the largest integer not exceeding the real
number $x$.
We remark that $R_{A_u}$ is ``close'' to $u$,
in the sense that $R_{A_u} \le u < R_{1+A_u}$.
In fact, the difference $u-R_{A_u}$ (known as the age of the renewal process)
is a tight family (over $u \ge 0$) of random variables.
In the above theorem, we just introduce another parameter, $t$,
and measure the difference between $tu$ and $R_{[tA_u]}$.
When $t=0$ or $1$, this difference is ``negligible''
with respect to any power of $u$. When $t$ is between $0$ and
$1$, then the difference is of the ``order of $\sqrt{u}$''
in the sense that when divided by $\sqrt{u}$ it converges to
a normal random variable. Jointly, over all $t \in [0,1]$,
we have convergence to a Brownian bridge, and this is what we
show next.

\proof
Consider, for $u > 0$,
\[
y_u(t) := \frac{R_{[tu]} - \mu tu}{\sqrt{u}}, \quad t \ge 0. 
\]
{From} Donsker's theorem \cite{BILL68}
for the random walk $\{R_n\}$
we have that $y_u \Rightarrow \sigma W$, where $W$ is a standard
Brownian motion.
Define also, for $u > 0$,
\[
\varphi_u(t) := \frac{tA_u}{u}.
\]
{From} the law of large numbers for the renewal process, $A_u/u \to
\mu^{-1}$, a.s., as $u \to \infty$. Hence, $\varphi_u$ converges
a.s.\ (and weakly) to the deterministic process $\{\mu^{-1} t\}$.
Since composition is a continuous function \cite{BILL68} we have that
\begin{equation}
\label{BMC}
\{(y_u \circ \phi_u) (t)\} \Rightarrow \{\sigma W_{\mu^{-1} t}\} 
\stackrel{d}{=} \{\mu^{-1/2}\sigma W_t\}.
\end{equation}
We also have
\[
(y_u \circ \varphi_u) (t)  =
\frac{R_{[tA_u]}-\mu t A_u}{\sqrt{u}},
\]
and so
\begin{align}
\eta_u(t) &= (y_u \circ \varphi_u) (t) 
+ \mu t \frac{A_u-\mu^{-1} u}{\sqrt{u}} \nonumber \\
&= (y_u \circ \varphi_u) (t) 
-t (y_u \circ \varphi_u)(1) -t \frac{u-R_{A_u}}{\sqrt{u}}.
\label{decomp}
\end{align}
Observe now that $\{u-R_{A_u}, u \ge 0\}$ is a tight family.
Indeed, from standard renewal theory (see, e.g., \cite{ASM03}),
if $R_1$ has a non-lattice distribution, then $u-R_{A_u}$ converges
weakly as $u\to\infty$. And if $R_1$ has a lattice distribution
with span $h$, then a similar convergence takes places for
$nh-R_{A_{nh}}$ as $n\to\infty$. Since, for all $u \ge 0$, 
$0 \le u-R_{A_u} \le ([u/h]+1)h-R_{A_{[u/h]}}$, the family
$\{u-R_{A_u}, u \ge 0\}$ is tight even in the lattice case.
Tightness implies that the last term of \eqref{decomp}
converges to $0$ in probability.
{From} the convergence stated in \eqref{BMC} and the 
decomposition \eqref{decomp}, we have that
\[
\{\eta_u(t)\}_{0 \le t \le 1} \Rightarrow 
\mu^{-1/2}\sigma \{W_t - t W_1\}_ {0 \le t \le 1}.
\]
It is well known \cite{REVYOR99}
that a standard Brownian bridge $W^0$ can be
represented as $W^0_t = W_t - t W_1$, and so the process above is 
the limit we were looking for.
\qed

\section{Extensions, discussion, and examples}
Here is a different version that, perhaps, makes Theorem \ref{BBthm} 
clearer:
Suppose that $M$ is a regenerative random measure on $[0, \infty)$.
That is, there is some renewal process with points $T_0 < T_1< T_2
< \cdots$ 
such that the random measures obtained by restricting $M$ onto 
$[T_n, T_{n+1})$, $n=0,1,2,\dots$, are i.i.d.
Suppose that 
\begin{gather*}
\mu:= E(T_2-T_1), \quad
\var (T_2-T_1) < \infty,
\\
\alpha := E M([T_1, T_2)) , \quad
0 < \var  M([T_1, T_2)) < \infty.
\end{gather*}
Define the random distribution function of $M$ by
\[
S(t) = M((0,t]), \quad u \ge 0.
\]
By the law of large numbers, $S(t)/t \to \mu^{-1} \alpha$, a.s.\ as
$t \to \infty$.
Consider the generalised inverse
\[
S^{-1}(u) :=\inf\{t \ge 0:~ S(t) > u \}, \quad u \ge 0.
\]
Then, in some naive sense,
$S^{-1}$ composed with $S$ is ``approximately'' the identity function,
but what can we say about the composition of $S^{-1}$ with a fraction
$t S$ of $S$ where $0< t < 1$?
The law of large numbers tells us that, almost surely,
\[
\frac{S(tS^{-1}(u))}{u} \xrightarrow[u \to \infty]{} t.
\]
An extension of the previous theorem quantifies the deviation:
\begin{theorem}
\label{extension}
As $u \to \infty$, the sequence of processes $\eta_u$ where
\[
\eta_u(t) := \frac{S(t S^{-1}(u)) - tu}{\sqrt{u}},
\quad 0 \le t \le 1,
\]
converges weakly to a Brownian bridge.
\end{theorem}
The proof of this is analogous to the previous one, so it is omitted.
Observe that the ``tying down'' of the Brownian motion occurs
naturally at $t=0$ and $t=1$.

The Brownian bridge has a scaling constant depending on
the parameters of the process $S$.

Note that the regenerative assumption is not crucial. All we need
is to have a process for which a Donsker theorem with a Brownian
limit holds. This is then translatable to a Brownian bridge limit.

If we interchange the roles of $S$ and $S^{-1}$ we still get a Brownian
bridge but with different constant.
For instance, interchanging the roles of $\{R_n\}$ and $\{A_u\}$ in
Theorem \ref{BBthm} we obtain that
\[
\eta'_n(t) := \frac{A(t R_n)-tn}{\sqrt{n}} , \quad 0 \le t \le 1,
\]
converges weakly, as $n \to \infty$, to $\kappa W^0$, where $W^0$
is a standard Brownian bridge and $\kappa = \sigma \mu^{-1}$.

\subsection{An interpretation}
To better understand the phenomenon, we cast the limit theorem 
as follows: We have a random function $S$, composed with scaling
functions
\[
\rho_t : x \mapsto t x
\]
and composed again with the inverse function $S^{-1}$
and we look at the asymptotic
behaviour
of the family of random functions
\begin{equation}
\label{hom}
S \comp \rho_t \comp S^{-1} - \rho_t, \quad 0 \le t \le 1,
\end{equation}
(or of $S^{-1} \comp \rho_t \comp S$),
as a function of the parameter $t$.
Thus, the time parameter of the Brownian bridge obtained in the limit plays
the role of a scaling factor.
When $t$ is $0$ or $1$, $S \comp \rho_t \comp S^{-1} - \rho_t$
is approximately zero (with respect to the normalising factor).
This raises the following three questions: 
\\
(i) How much ``one-dimensional'' is this phenomenon?
\\
(ii) Can we replace the family $\rho_t$ by a more general homotopy?
\\
(iii) Are different kind of bridges possible to obtain?
\\
With respect to the latter question, we could start with a regenerative 
process with finite mean but infinite variance, one that belongs to
the domain of attraction of, say, a self-similar L\'evy process.

\subsection{Four examples}

\paragraph{EXAMPLE 1}
The first is a simple example involving a standard Brownian motion $W$. 
Let $X$ denote the (strong) Markov process
\begin{equation}
\label{reflW}
X_t = (W_t-t)-\min_{0 \le s \le t}(W_s-s), \quad t \ge 0,
\end{equation}
which is the reflection of the drifted Brownian motion $\{W_t-t\}$.
This process in natural in many areas of applied probability, e.g.\
in the diffusion approximation of a queue.
We have $X_0=0$, $X_t \ge 0$. The {\em Brownian area process}
\begin{equation}
\label{areaX}
S(t) = \int_0^t X_r dr
\end{equation}
is non-decreasing. 
Fix some $u \ge 0$ and $t \in [0,1]$.
By continuity, there is a unique point between $0$ and $u$
that splits the area $S(u)$ into two parts with ratio $t:(1-t)$.
Call this point $H_u(t)$.
Specifically,
\[
H_u(t) :=
\min\left\{v \ge 0:~ t \int_0^v X_r dr = (1-t) \int_v^u X_r dr\right\},
\quad 0 \le t \le 1.
\]
We then claim that
\[
\eta_u(t) := \frac{H_u(t) - tu}{\sqrt{u}}, \quad 0 \le t \le 1,
\]
converges weakly to a Brownian bridge
as $u \to \infty$.
To see this, observe that
\[
S^{-1}(x) = \min\{v \ge 0:~ S(v) = x\},
\]
and hence
\begin{align*}
S^{-1}(t S(u))
&= \min\{v \ge 0:~ S(v) = t S(u)
\\
&= \min\{v \ge 0:~ S(v) = t (S(v)+S(u)-S(v))\}
\\
&= \min\{v \ge 0:~ (1-t) S(v) = t (S(u)-S(v))\} = H_u(1-t).
\end{align*}
Apply Theorem \ref{extension} to get the result.
(Notice that $\eta_u(1-t)$ also converges to a Brownian bridge.)

\paragraph{EXAMPLE 2}
Same as Example 1, but with $W$ being a {\em zero-mean} L\'evy process.
The Brownian bridge in Example 1 was obtained not from the fact that
$W$ was Brownian, but from the regenerative structure of $S$.
It is this that allows us to replace $W$ by a more general, say a L\'evy
process, as long as we maintain the finite variance assumptions.
The latter hold once we add a strictly negative drift to
a zero-mean L\'evy process $W$, reflect it, precisely as in 
\eqref{reflW}, and integrate just as in \eqref{areaX}.
Whereas $W$ may be discontinuous, $S$ is continuous and
the conclusion remains the same.

\paragraph{EXAMPLE 3}
The third example is an application of the above in proving a
limit theorem for a random digraph.
We consider a random directed graph $G_n=(V_n, E_n)$ on the set of vertices 
$V_n:=\{1, \ldots, n\}$
by letting the set of edges $E_n$ 
contain the pair $(i,j)$, $i < j$, with probability $p$, 
independently from pair to pair.
This is a directed version of the (nowadays) so-called Erd\H{o}s-R\'enyi graph.

A path starting in $i$ and ending in $j$ is a sequence of vertices
$i_0=i, i_1, \ldots, i_n=j$ such that
$(i,i_1), \ldots (i_{n-1},j)$ are edges.
Amongst all paths in $G_n$ there is one with maximum length; this length is
denoted by $L_n$.
Amongst all paths in $G_n$ that end at a vertex $j \in V_n$
there is one with maximum length; this length is called {\em weight} of 
vertex $j$.
We keep track of vertices with a specific weight and let
$S_n(\ell)$ be the number of vertices with weights {\em at least} $\ell$.
(Here $\ell$ ranges between $0$ and $L_n$.)
So, for example, $S_n(0)$ is the number of vertices in $V_n$ that are
endpoints of no edge in $E_n$, and $S_n(L_n)$ is the number of paths of maximal
length in $G_n$.

\begin{theorem}
\[
\frac{S_n([t L_n]) - tn}{\sqrt{n}}, \quad 0 \le t \le 1,
\]
converges, as $n \to \infty$, weakly to a Brownian bridge.
\end{theorem}
The proof of this theorem can be found in \cite[p.\ 453]{FK03}.

\paragraph{EXAMPLE 4}
Here is an illustration, of the kind of phenomenon described around \eqref{hom},
in Stochastic Geometry.
We consider a Poisson point process\footnote{More general point 
processes can be allowed here.} 
$N$ in $\R^d$ with intensity, say, $1$;
that is, $N$ is a random discrete subset of $\R^d$ such that
the cardinalities of $N\cap B_1, \ldots, N \cap B_n$ are independent
random variables whenever $B_1, \ldots, B_n$ are disjoint Borel sets,
for any $n \in \N$, and the expectation of the cardinality of $N \cap B$
equals the Lebesgue measure of $B$.
For each $x$ in $\R^d$ we let $\pi(x)$ be the point of $N$ closest to $x$
(there is a.s.\ a unique such point).
For each point $z$ of $N$, we  let $\sigma(z)$ be 
the {\em Voronoi cell} \cite{SW00,SKM87} associated to $z$:
\[
\sigma(z) := \{x \in \R^d:~ ||x-z|| \le ||x-z'|| \text{ for all points $z'$
of $N$}\},
\]
where $||\cdot||$ is the Euclidean norm on $\R^d$.
The {\em Voronoi tessellation} of $\R^d$ is the
the tiling of $\R^d$ by the Voronoi cells.
If $z$ is not a point of $N$ we define $\sigma(z)$ to be the Voronoi cell
containing $z$ (again this cell is a.s.\ unique).
The distance of a closed set $A \subset \R^d$ from a point $x \in \R^d$ is
\[
\operatorname{dist}(A, x) = \inf\{||x-y||:~ y \in A\}.
\]
Consider now the process
\[
D(t, x) := \operatorname{dist}(\sigma(t \pi(x)), t x),
\]
where $t \in [0,1]$ and $x \in \R^d$.
The claim is that
\[
||x||^{-1/2}~D(\cdot,x) \Rightarrow |W^0|, \quad \text{ as } ||x|| \to \infty,
\]
$|W^0|$ being the absolute value of a Brownian bridge.

\vspace*{5mm}
{\sc School of Mathematical and Computer Sciences\\
Maxwell Institute for Mathematical Sciences\\
Heriot-Watt University, Edinburgh
EH14 4AS, U.K.}\\
{\tt foss@ma.hw.ac.uk},
{\tt takis@ma.hw.ac.uk}

\end{document}